\documentstyle{amlts}
\begin{document}
\annalsline{158}{2003}
\received{July 23, 2001}
\revised{October 21, 2002}
\startingpage{661}
\def\bye{\end{document}}
 \font\tenrm=cmr10
\def\ritem#1{\item[{\rm #1}]}
\input amssym.def
\input amssym.tex

 %\input BoxedEps.Tex %for mac
%\SetTexturesEPSFSpecial %for mac
\input boxedeps.tex % unix
\SetepsfEPSFSpecial % unix
\HideDisplacementBoxes
\def\figin#1#2{
$$
 {\BoxedEPSF{#1 scaled
#2}%
}%
$$
\noindent}
\def\scr#1{{\scriptscriptstyle #1}}
\def\eqref#1{(\ref{#1})}
%--------------- Author macros ---------------
%for Bbb in amstex
\catcode`\@=11
\font\twelvemsb=msbm10 scaled 1100
\font\tenmsb=msbm10
%\font\ninemsb=msbm7 scaled 1100%msbm9
\font\ninemsb=msbm10 scaled 800
\newfam\msbfam
\textfont\msbfam=\twelvemsb  \scriptfont\msbfam=\ninemsb
  \scriptscriptfont\msbfam=\ninemsb
\def\msb@{\hexnumber@\msbfam}
\def\Bbb{\relax\ifmmode\let\next\Bbb@\else
 \def\next{\errmessage{Use \string\Bbb\space only in math
mode}}\fi\next}
\def\Bbb@#1{{\Bbb@@{#1}}}
\def\Bbb@@#1{\fam\msbfam#1}
\catcode`\@=12

 \catcode`\@=11
\font\twelveeuf=eufm10 scaled 1100
\font\teneuf=eufm10
\font\nineeuf=eufm7 scaled 1100%eufm9
\newfam\euffam
\textfont\euffam=\twelveeuf  \scriptfont\euffam=\teneuf
  \scriptscriptfont\euffam=\nineeuf
\def\euf@{\hexnumber@\euffam}
\def\frak{\relax\ifmmode\let\next\frak@\else
 \def\next{\errmessage{Use \string\frak\space only in math
mode}}\fi\next}
\def\frak@#1{{\frak@@{#1}}}
\def\frak@@#1{\fam\euffam#1}
\catcode`\@=12
%-------------- Author entries --------------------
\renewcommand{\leq}{\leqslant}
\renewcommand{\geq}{\geqslant}
\newcommand{\maps}[3]{{#1}\!:\!{#2}\rightarrow{#3}}
\newcommand{\floor}[1]{{\lfloor{#1}\rfloor}}

\newcommand{\Hdim}{\mathop{\rm H.dim}}
\newcommand{\NE}{\mathop{\rm NE}}
\newcommand{\Spec}{\mathop{\rm Spec}}
\newcommand{\dens}{\mathop{\rm dens}}
\newcommand{\area}{\mathop{\rm area}}
\newcommand{\dist}{\mathop{\rm dist}}
\newcommand{\myliminf}{\mathop{\rm varliminf}}
\newcommand{\srliminf}{\mathop{\rm liminf}}
\newcommand{\bbZ}{{\Bbb Z}}

\title{Hausdorff dimension of the set\\ of nonergodic directions}
\shorttitle{Hausdorff dimension of nonergodic directions}  

 \author{Yitwah Cheung}
 \institutions{Rice University, Houston, TX\\
{\eightpoint {\it E-mail address\/}: michael@math.rice.edu}}

\vglue-16pt
 \centerline {\tenrm (with an Appendix by {\tensc M. Boshernitzan})}

\vglue36pt \centerline{\bf Abstract}\vglue12pt

It is known that nonergodic directions in a rational billiard form 
a subset of the unit circle with Hausdorff dimension at most $1/2$.  
Explicit examples realizing the dimension $1/2$ are constructed 
using Diophantine numbers and continued fractions.  A lower 
estimate on the number of primitive lattice points in certain 
subsets of the plane is used in the construction.

\section{Introduction}\label{S:Intro}

Consider the billiard in a polygon $Q$.  A fundamental result \cite{KMS} 
implies that a typical trajectory with typical initial direction will be 
equidistributed provided the angles of $Q$ are rational multiples of $\pi$.  
More precisely, there is a flat surface $X$ associated to the polygon 
such that each direction $\theta\in S^1$ determines an area-preserving 
flow on $X$; the assertion is that the set ${\rm NE}(Q)$ of parameters $\theta$ 
for which the associated flow is not ergodic has measure zero.  The 
statement holds more generally for the class of \emph{rational billiards} 
in which the (abstract) polygon is assumed to have the property that 
the subgroup of $O(2)$ generated by the linear parts of the reflections 
in the sides is finite.  For a recent survey of rational billiards, 
see \cite{MT}.

Let $Q_\lambda, \lambda\in(0,1)$, be the polygon described informally 
as a $2$-by-$1$ rectangle with an interior wall extending orthogonally 
from the midpoint of a longer side so that its distance from the opposite 
side is exactly $\lambda$ (see Figure~1). 
We are interested in the Hausdorff dimension of the set $\NE(Q_\lambda)$.  
Recall that $\lambda$ is \emph{Diophantine} if the inequality 
$$\left|\lambda-\frac{p}{q}\right|\leq\frac{1}{|q|^e}$$ has 
(at most) finitely many integer solutions for some exponent $e>0$. 

\centerline{\BoxedEPSF{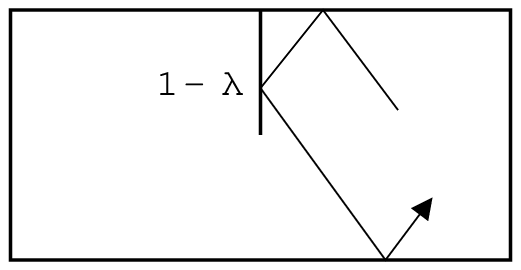 scaled 1000}}
 
\centerline{Figure 1. The billiard in $Q_\lambda$.}
\vglue28pt

\specialnumber{1}\proclaim{Theorem} \label{Thm1}
If $\lambda$ is Diophantine{\rm ,} then $\Hdim \NE(Q_\lambda)=1/2$.  
\endproclaim

\vglue-12pt

In fact, Masur has shown that for any rational billiard the set of 
nonergodic directions has Hausdorff dimension at most $1/2$ \cite{Ma}.  
This upperbound is sharp, as Theorem~\ref{Thm1} shows.  It should be 
pointed out that the theorems in \cite{KMS} and \cite{Ma} are stated 
for holomorphic quadratic differentials on compact Riemann surfaces.  
The flat structure on the surface associated to a rational billiard 
is a special case, namely the square of a holomorphic $1$-form.  

The ergodic theory of the billiards $Q_\lambda$ was first studied 
by Veech \cite{V1} in the context of $\bbZ_2$ skew products of 
irrational rotations.  Veech proved the slope of the initial 
direction $\theta$ has bounded partial quotients if and only if the corresponding 
flow is (uniquely) ergodic for all $\lambda$.  
On the other hand, if $\theta$ has unbounded partial quotients, then 
there exists an uncountable set $K(\theta)$ of $\lambda$ for which 
the flow is not ergodic.  In this way, Veech showed that minimality 
does not imply (unique) ergodicity for these $\bbZ_2$ skew products.  
(The first examples of minimal but uniquely ergodic systems had been 
constructed by Furstenberg in \cite{Fu}.)  Our approach is dual to 
that of Veech in the sense that we fix $\lambda$ and study the set 
of paramaters $\theta\in\NE(Q_\lambda)$.  

The billiards $Q_\lambda$ were first introduced by Masur and Smillie 
to give a geometric representation of the $\bbZ_2$ skew products 
studied by Veech.  It follows from \cite{V1} that $\NE(Q_\lambda)$ 
is countable if $\lambda$ is rational.  A proof of the converse can 
be found in the survey article \cite[Thm.\ 3.2]{MT}.  Boshernitzan has given 
a short argument showing $\Hdim\NE(Q_\lambda)=0$ for a residual 
(hence, uncountable) set of $\lambda$.  (His argument is presented 
in the appendix to this paper.)  Theorem~\ref{Thm1} implies any 
such $\lambda$ is a Liouville number.  As is well-known, the set 
of Liouville numbers has measure zero (in fact, Hausdorff dimension 
zero).  We remark that by Roth's theorem every algebraic integer 
satisfies the hypothesis of Theorem~\ref{Thm1}.  

Some generalizations of Theorem~\ref{Thm1} are mentioned in Section~\ref{S:SumX}.  
For the class of Veech billiards (see \cite{V2}) the set of nonergodic 
directions is countable.  It would be interesting to know if there are 
(number-theoretic) conditions on a general rational billiard $Q$ which 
imply that the Hausdorff dimension of $\NE(Q)=1/2$.

Theorem~\ref{Thm1} can be reduced to a purely number-theoretic statement.

\proclaimtitle{Summable cross products condition}
\proclaim{Lemma}\label{L:SumX}
Suppose $(w_j)$ is a sequence of vectors of the form 
$(\lambda+m_j,n_j)${\rm ,} where $m_j,n_j\in2{\Bbb Z}$ and $n_j\neq0${\rm ,} and 
assume that the Euclidean lengths $|w_j|$ are increasing.  The condition 
\begin{equation}\label{E:SumX}
	\sum \left| w_j \times w_{j+1} \right| < \infty, 
\end{equation}
implies that $\theta_j = w_j/|w_j|$ converges to some 
$\theta\in\NE(Q^t_\lambda)$ as $j\to\infty$.  
{\rm (}\/Here{\rm ,} $Q^t_\lambda$ is the billiard table obtained 
by reflecting $Q_\lambda$ in a line of slope $-1$.{\rm )}   
\endproclaim

\specialnumber{2}\proclaim{Theorem} \label{Thm2}
Let $K(\lambda)$ be the set of nonergodic directions 
that can be obtain using Lemma~{\rm \ref{L:SumX}.}
If $\lambda$ is Diophantine{\rm ,} then $\Hdim K(\lambda)=1/2$.  
\endproclaim

\demo{Proof of Theorem~{\rm \ref{Thm1}}}
Theorem~\ref{Thm2} implies $\Hdim \NE(Q_\lambda)=\Hdim 
\NE(Q^t_\lambda)\break\geq1/2$.  Together with Masur's upperbound, 
this gives Theorem~\ref{Thm1}.  
\enddemo

{\it Density of primitive lattice points}.
The main obstacle in our approach to finding lowerbounds on Hausdorff 
dimension is the absence of \emph{primitive} lattice points in certain 
regions of the plane.  More precisely, let $\Sigma=\Sigma(\alpha,R,Q)$ 
denote the parallelogram (Figure~2) 
$$\Sigma := \left\{ (x,y)\in{\Bbb R  }^2 
	: |y\alpha-x|\leq1/Q, \; R \leq y \leq 2R \right\}$$
and define 
$$\dens(\Sigma):=\frac{^\#\{(p,q)\in\Sigma : \gcd(p,q)=1\}}{\area(\Sigma)}.$$
\figin{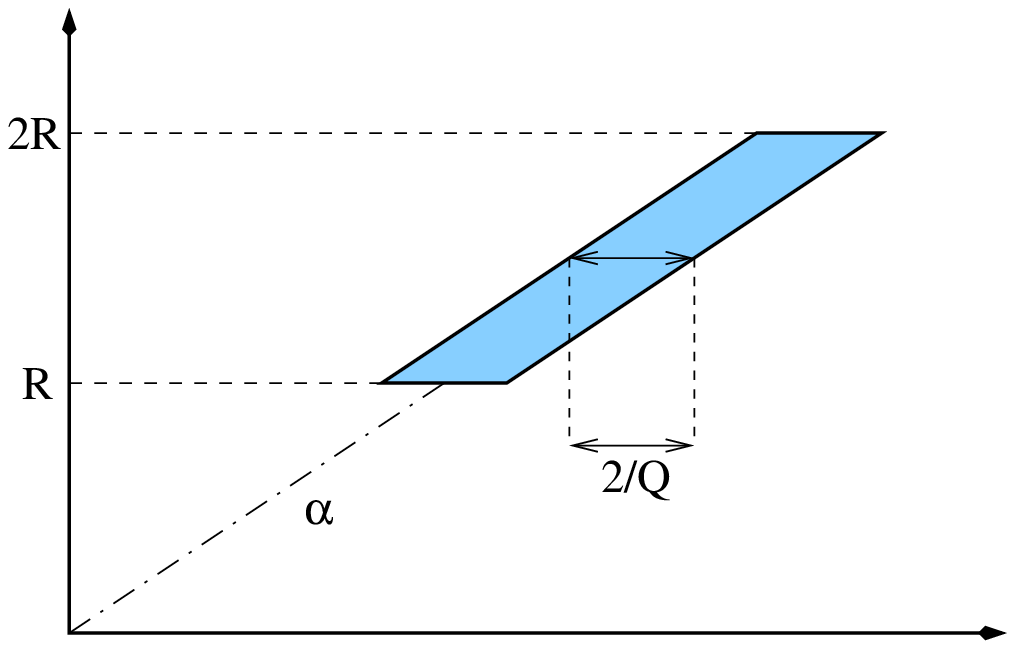}{700}

\centerline{Figure 2. The parallelogram $\Sigma(\alpha,R,Q)$.}
\vglue12pt
 
The proof of Theorem~\ref{Thm2} relies on the following fact:  

\specialnumber{3}\proclaim{Theorem} \label{Thm3} 
Let $\Spec(\alpha)$ be the sequence of heights formed by 
the convergents of $\alpha$.  There exist constants $A_0$ 
and $\rho_0>0$ such that whenever $\area(\Sigma)\geq A_0$ 
\vglue9pt
\hfill ${\displaystyle \Spec(\alpha)\cap[Q,R] \neq \emptyset 
  \quad\; \Rightarrow \quad \dens(\Sigma)\geq\rho_0.}$\hfill
\endproclaim

{\it Remark.} 
It can be shown $\dens(\Sigma)=0$ if $\alpha$ does not 
have any convergent whose height is between $Q/4$ and $8R$.  
Thus, $\area(\Sigma)\gg1$ alone cannot imply the existence 
of a primitive lattice point in $\Sigma$.  For example, the 
implication $$|\alpha|<\frac{1}{2R}\left(1-\frac{1}{Q}\right) 
\quad \Rightarrow \quad \dens(\Sigma)=0$$ is easy to verify 
and remains valid even if $|\cdot|$ is replaced by the distance 
to the nearest integer (because arithmetic density is preserved 
under $ (\hskip-4pt\begin{array}{cc} \scr{1}&\hskip-6pt\scr{n}\\[-8pt]
\scr{0}&\hskip-6pt\scr{1}\end{array}\hskip-4pt )$).

\demo{Outline}
Theorem~\ref{Thm2} is proved by showing that $K(\lambda)$ contains 
a Cantor set whose Hausdorff dimension may be chosen close to $1/2$ 
when $\lambda$ is Diophantine.  The construction of this Cantor set 
is based on Lemma~\ref{L:SumX} and is presented in Section~\ref{S:SumX}.  
The proof of Theorem~\ref{Thm2} is completed in Section~\ref{S:Diophantine} 
if we assume the statement of Theorem~\ref{Thm3}, whose proof is deferred 
to Section~\ref{S:Density}.  
\enddemo

{\it Acknowledgments.}
This research was partially supported by the National Science Foundation 
and the Clay Mathematics Institute.  The author would also like to 
thank his thesis advisor Howard Masur for his excellent guidance.

\section{Cantor set of nonergodic directions}\label{S:SumX}

We begin with the proof of Lemma~\ref{L:SumX}, which is the recipe 
for the construction of a Cantor set $E(\lambda)\subset K(\lambda)$.  
We then show that the Hausdorff dimension of $E(\lambda)$ can be 
chosen arbitrarily close to $1/2$ if the arithmetic density of 
the parallelograms $\Sigma(\alpha,R,Q)$ can be bounded uniformly 
away from zero.  

\demo{{\rm 2.1.}  Partition determined by a slit}
The flat surface associated to $Q_\lambda$ is shown in Figure~3.  It will be slightly more convenient to
work with the  reflected table $Q^t_\lambda$.  Let $X_\lambda$ be the flat surface 
associated to $Q^t_\lambda$.  The proof of Lemma~\ref{L:SumX} is 
based on the following observation:  
$$X_\lambda \hbox{ is a branched double cover of the 
	square torus } T={\Bbb R}^2/{\Bbb Z}^2.$$
More specifically, let $w_0\subset T$ denote the projection of the 
interval $[0,\lambda]$ contained in the $x$-axis.  $X_\lambda$ may 
be realized (up to a scale factor of $2$) by gluing two copies of 
the slit torus $T\setminus w_0$ along their boundaries so that the 
upper edge of the slit in one copy is attached to the lower edge 
of the slit in the other, and similarly for the remaining edges.  
The induced map $\maps{\pi}{X_\lambda}{T}$ is the branched double 
cover obtained by making a cut along the slit $w_0$.  \enddemo

\proclaimtitle{Slit directions are nonergodic}
\proclaim{Lemma}\label{L:Slit}
A vector of the form\break $(\lambda+m,n)$ with $m,n\in2{\Bbb Z}$ 
and $n\neq0$ determines a nonergodic direction in~$Q^t_\lambda$.  
\endproclaim

\centerline{\BoxedEPSF{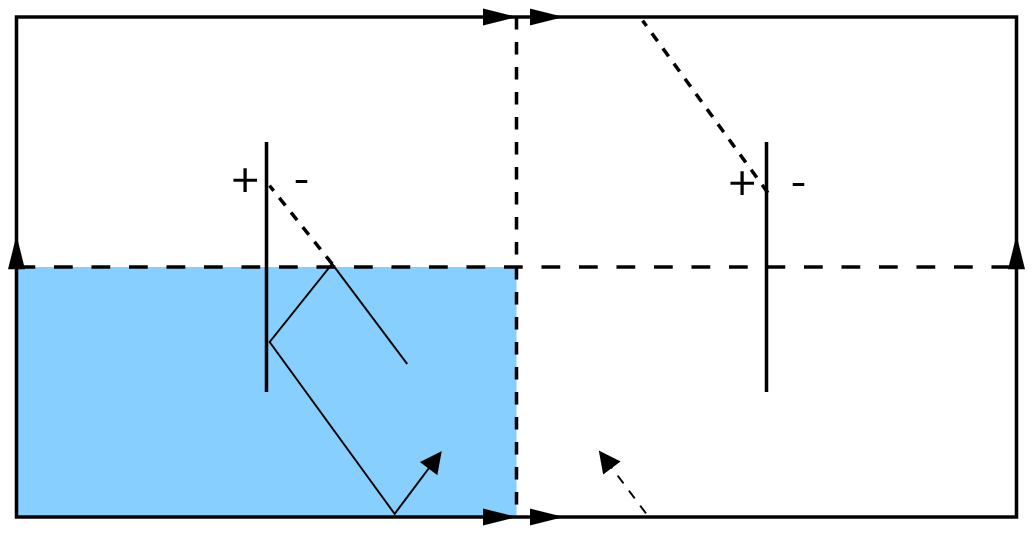 scaled 900}}
\vglue6pt
\centerline{Figure 3. Unfolded billiard trajectory.}

\vglue12pt
{\it Proof}.
A vector of the given form determines a slit $w$ in $T$ that is 
homologous to $w_0$ (mod 2).  (We assume $\lambda$ is irrational, 
for the statement of the lemma is easily seen to hold otherwise.)  
If $\maps{\pi'}{X'}{T}$ is the branched double cover obtained by 
making a cut along $w$, then there is a biholomorphic isomorphism 
$\maps{h}{X_\lambda}{X'}$ such that $\pi=\pi'\circ h$.   It follows 
that $\pi^{-1}(w)$ partitions $X_\lambda$ into a pair of slit tori 
with equal area, and that this partition is invariant under the 
flow in the direction of the slit.  Hence, the vector (after 
normalization) determines a nonergodic direction in $Q^t_\lambda$.  
\hfill\qed 

\demo{Proof  of Lemma~{\rm \ref{L:SumX}}}
It is easy to see from (\ref{E:SumX}) that the directions $\theta_j$ 
form a Cauchy sequence.  The corresponding partitions of $X_\lambda$ 
also converge in a measure-theoretic sense: the symmetric difference of 
consecutive partitions is a union of parallelograms whose total area is 
bounded by the corresponding term in (\ref{E:SumX}); summability implies 
the existence of a limit partition.  Invariance of the limit partition 
under the flow in the direction of $\theta$ will follow by showing that 
$h_j$, the component of $w_j$ perpendicular to $\theta$, tends to zero 
as $j\to\infty$ (\cite[Th.~2.1]{MS}).  To see this, observe that 
the area of the right triangle formed by $w_j$ and $\theta$ is roughly 
$h_j$ times the Euclidean length of $w_j$; it is bounded by the tail 
in (\ref{E:SumX}) and therefore tends to zero.  
(We have implicitly assumed that $\lambda$ is irrational.  For rational 
$\lambda$ the lemma still holds because a nonzero term in (\ref{E:SumX}) 
must be at least the reciprocal of the height.)  
\enddemo

{\it Remark.}
A vector of the form $(\lambda+m,n)$ with $m,n\in g{\Bbb Z}$ and 
$n\neq0$ determines a partition of the branched $g$-cyclic cover of 
$T$ into $g$ slit tori of equal area.  From this, it is not hard to 
show that the conclusion of Theorem~\ref{Thm1} holds in genus $g\geq2$.  
Gutkin has pointed out other higher genus examples obtained by 
considering branched double covers along multiple parallel slits.  
Further examples are possible by observing that the proof of 
Theorem~\ref{Thm2} depends only on a Diophantine condition 
on the vector $w_0=(\lambda,0)$. (See \S\ref{S:Diophantine}.)

\demo{{\rm 2.2.}  Definition of $E(\lambda)$}
Our goal is to find sequences that satisfy condition (\ref{E:SumX}) 
and intuitively, the more we find, the larger the dimension.  However, 
in order to facilitate the computation of Hausdorff dimension, 
we shall restrict our attention to sequences whose Euclidean 
lengths grow at some fixed rate.  

We shall realize $E(\lambda)$ as a decreasing intersection of compact 
sets $E_j$, each of which is a disjoint union of closed intervals.  
Let $V$ denote the set of vectors that satisfy the hypothesis of 
Lemma~\ref{L:Slit}.  Henceforth, by a \emph{slit} we mean a vector 
$w\in V$ whose \emph{length} is given by $L:=|n|$ and \emph{slope} 
by $\alpha:=(\lambda+m)/n$.  Note that the following version of the 
cross product formula holds: $|w\times w'|=LL'\Delta$, where $\Delta$ 
is the distance between the slopes.  Fix a parameter $\delta>0$.  
\enddemo

\advance\theoremcount by 1
{\it Definition}   2.2  (Children of a slit).
Let $w$ be a slit of length $L$ and slope~$\alpha$.  
A slit $w'$ is said to be a \emph{child} of $w$ if 
  \begin{itemize}
    \item[(i)] $w'=w+2(p,q)$ for some relatively prime integers $p$ and $q$ 
    \item[(ii)] $|q\alpha-p| \leq 1/L\log L$ and 
		$q\in[L^{1+\delta},2L^{1+\delta}]$.
  \end{itemize}
 
\vglue4pt
\proclaimtitle{Chains have nonergodic limit}
\proclaim{Lemma}\label{L:Chain}
The direction of $w_j$ converges to a point in $K(\lambda)$ as 
$j\to\infty$ provided $w_{j+1}$ is a child of $w_j$ for every~$j$. 
\endproclaim

\vglue-16pt
{\it Proof}.
The inequality in (ii) (equivalent to $|w\times w'|\leq1/\log L$) 
implies that the directions of the slits are close to one another.  
Hence, their Euclidean lengths are increasing since the length of 
a child is approximately $L^{1+\delta}$.  The sum in (\ref{E:SumX}) 
is dominated by a geometric series of ratio $1/(1+\delta)$.  
\hfill\qed\vglue9pt

Choose a slit $w_0$ and call it the slit of level $0$.  The slits 
of level $j+1$ are defined to be children of slits of level $j$.  
Let $V':=\cup V'_j$ where $V'_j$ denotes the collection of slits 
that belong to level $j$.  Associate to each $w\in V'$ the smallest 
closed interval containing all the limits obtainable by applying 
Lemma~\ref{L:Chain} to a sequence beginning with $w$.  Define 
$E(\lambda):=\cap E_j$ where $E_j$ is the union of the intervals 
associated to slits in $V'_j$.  It is easily seen that the diameters 
of intervals in $E_j$ tend to zero as $j\to\infty$.  Hence, every 
point of $E$ arises as the limit obtained by an application of 
Lemma~\ref{L:Chain}.  Therefore, $E(\lambda)\subset K(\lambda)$.  

\vglue6pt  2.3. {\it Computation of Hausdorff dimension}.
We first give a heuristic calculation which shows that the Hausdorff 
dimension of $K(\lambda)$ is at most $1/2$.  (This fact is not used 
in the proof of Theorem~\ref{Thm1}.)  We then show rigorously that 
the Hausdorff dimension of $E(\lambda)$ is at least $1/2$ under 
a critical assumption: \emph{each slit in $V'$ has enough children}.  

Recall the construction of the Cantor middle-third set.  At each 
stage of the induction, intervals of length $\Delta$ are replaced 
with $m=2$ equally spaced subintervals of common length $\Delta'$.  
In this case, the Hausdorff dimension is exactly $\log2/\log 3$, 
or $\log m/\log(1/\varepsilon)$ where $\varepsilon:=\Delta'/\Delta=1/3$.  

For $K(\lambda)$ it is enough to consider sequences for which every 
term in (\ref{E:SumX}) is bounded above.  Associated to each slit 
of length $L$ is an interval of length $\Delta=1/L^2$.  The number 
of slits of length approximately $L'$ is at most $m=L'/L$.  Their 
intervals have approximate length $\Delta'=1/(L')^2$.  Therefore, $$\Hdim 
K(\lambda) \leq \frac{\log m}{\log (\Delta/\Delta')} = \frac{1}{2}.$$

To get a lowerbound on the Hausdorff dimension of $E(\lambda)$ we 
need to show there are lots of children and wide gaps between them.  
The number of children is exactly $2L^\delta/\log L$ times the 
arithmetic density of the parallelogram $\Sigma(\alpha,R,Q)$ 
where $R=L^{1+\delta}$ and $Q=L\log L$.  

\proclaimtitle{Slopes of children are far apart}
\proclaim{Lemma}\label{L:Gaps}
The slopes of any two children of a slit with length $L$ are 
separated by a distance of at least $O(1/L^{2+2\delta})$.  
\endproclaim

\vglue-12pt
{\it Proof}.
Let $w$ be a slit of length $L$.  A child $w'$ has the form $w'=w+2v$ 
for some $v=(p,q)$.  If $w''=w+2v'$ is another child with $v'=(p',q')$, 
then $v'\neq v$.  Since both pairs are relatively prime, 
$|p/q-p'/q'|\geq1/qq'\geq1/4L^{2+2\delta}$.  The lemma follows 
by observing that the slope of $w'$ satisfies 
\vglue12pt
\hfill ${\displaystyle \left|\alpha'-\frac{p}{q}\right| = 
	\frac{|w'\times v|}{L'q} = \frac{|w\times v|}{(L+2q)q} \leq 
	\frac{L|q\alpha-p|}{2q^2} \leq \frac{1}{2L^{2+2\delta}\log L}.}$ 
\hfill\qed

\phantom{odd}
\proclaimtitle{Enough children implies dimension $1/2$}
\proclaim{Proposition}\label{P:Hdim}
Suppose there exists $c_1>0$ such that every slit in $V'$ 
has at least $c_1L^\delta/\log L$ children in $V'${\rm ,} where $L$ 
denotes the length of the slit.  Then $\Hdim K(\lambda)=1/2$.  
\endproclaim

{\it Proof}.
The length of a slit in $V'_j$ is roughly $L_j=L_0^{(1+\delta)^j}$, 
where $L_0$ denotes the length of the initial slit $w_0$.  
The number of children is at least $m_j = c_1L_j^\delta/\log L_j$ 
and their slopes are at least $\varepsilon_j=1/4L_j^{2+2\delta}$ apart.  
It follows by well-known estimates for computing Hausdorff dimension 
(we use \cite[Ex.~4.6]{Fa}) that 
$$\Hdim E(\lambda) \geq \srliminf_{j\to\infty} 
	\frac{\log(m_0\cdots m_{j-1})}{-\log(m_j\varepsilon_j)} 
    = \srliminf_{j\to\infty} 
	\frac{\sum_{i=0}^{j-1}\delta\log L_i}{(2+\delta)\log L_j} 
    = \frac{1}{2+\delta}.$$  
Together with the upperbound on $K(\lambda)$, this proves the lemma.  
\hfill\qed\vglue12pt

{\it Remark.}
Theorem~\ref{Thm3} allows us to determine when a slit has 
enough children.  It should by pointed out that Diophantine 
$\lambda$ does {\it not} imply every slit will have enough 
children.  We shall show that Proposition~\ref{P:Hdim} holds 
if $V'$ is replaced by a suitable subset.  (By the remark 
following Theorem~\ref{Thm3} one can easily show there are slits 
that do not have any children and whose directions form a dense set.)

\section{Diophantine condition}\label{S:Diophantine}

Let $w_0$ be the initial slit in the definition of $E(\lambda)$.  
The hypothesis that $\lambda$ is Diophantine implies there are 
constants $e_0>0$ and $c_0>0$ such that $$||w_0\times v|| = 
\min_{n\in{\Bbb Z}} |w_0\times v-n| \geq \frac{c_0}{|v|^{e_0}} 
\quad\; \hbox{for all $v\in{\Bbb Z}^2$, $v\neq0$.}$$
Fix a real number $N$ so that $e_0<N\delta$.  We assume the length of $w_0$ 
is at least some predetermined value $L_0=L_0(\lambda,\delta,N,e_0,c_0)$.  

\demo{Definition {\rm 3.1  (Normal slits)}} A slit of length $L$ and slope $\alpha$ is 
said to be {\em normal} if for every real number $n, 1\leq n\leq N+1$, 
\vglue12pt
\hfill ${\displaystyle\Spec(\alpha)\cap[e^{n\delta}L\log L,L^{1+n\delta}]\neq\emptyset.}$\hfill  
\enddemo

Let $V''$ be the subset of $V'$ formed by normal slits of length $\geq L_0$.  

\advance\theoremcount by 1

\proclaimtitle{Normal slits have enough children}
\proclaim{Proposition}\label{P:Normal} \hskip-8pt
There exists $c_1>0$ such that every slit in $V''$ has 
at least $c_1L^\delta/\log L$ children in $V''$.  
\endproclaim

To complete the proof of Theorem~\ref{Thm2} we also need 

\proclaimtitle{Normal slits exist}
\proclaim{Lemma}\label{L:Exist}
Arbitrarily long normal slits exist.  
\endproclaim

{\it Proof  of Theorem~{\rm \ref{Thm2}}  assuming 
 Lemma~{\rm \ref{L:Exist}} and  Proposition~{\rm \ref{P:Normal}}}.
We may choose the initial slit $w_0$ to lie in $V''$, which is nonempty 
by the lemma.  The calculation in the proof of Proposition~\ref{P:Hdim} 
applies to a subset of $E(\lambda)$ to give the same conclusion; 
in other words, the proposition implies $\Hdim K(\lambda)=1/2$.  \phantom{rain}\hfill\qed 
\vglue6pt

\advance\eqcount by 1
We recall two classical results from the theory of continued fractions.  
The $k^{\rm th}$ convergent $p_k/q_k$ of a real number $\alpha$ is 
a (reduced) fraction such that %which approximates $\alpha$ to within 
\begin{equation}\label{CFA}
	\frac{1}{q_k(q_{k+1}+q_k)} \leq
	\left|\alpha-\frac{p_k}{q_k}\right|
	\leq \frac{1}{q_kq_{k+1}}
\end{equation}
and satisfies the recurrence relation $q_{k+1}=a_{k+1}q_k+q_{k-1}$ 
(similarly for $p_k$), where $a_k$ is the $k^{\rm th}$ partial quotient.  
A partial converse is that if $p$ and $q>0$ are integers satisfying 
\begin{equation}\label{SCC}
\left|\alpha-\frac{p}{q}\right| \leq \frac{1}{2q^2}
\end{equation}%apart from a trivial exception
then $p/q$ is a convergent of $\alpha$, although it need not be reduced.   

\demo{{\rm 3.1.}  Existence of normal slits}
\vglue6pt
\advance\theoremcount by 1
{\it Definition} 3.4.
A slit of length $L$ and slope $\alpha$ is said to be $n$-{\it good} if 
$$\Spec(\alpha)\cap[e^{n\delta}L\log L,L^{1+\delta}]\neq\emptyset.$$  
 
\proclaim{Lemma}
\label{L:Good=>Normal}%%%%%%%%%%%%%%%%%%%%%%%%%%%%%%%%%%%%%%%%%%%%%%%%%
A sufficiently long $N$\/{\rm -}\/good slit is normal.  
\endproclaim

{\it Proof}.
An $(N+1)$-good slit is normal by definition, so it suffices to 
consider the case of an $N$-good slit that is not $(N+1)$-good.  
Suppose $w$ is such a slit, with length $L$ and slope $\alpha$.  
Let $q_k$ be the largest height in $\Spec(\alpha)\cap[1,L^{1+\delta}]$ 
so that $q_k=e^{n_1\delta}L\log L$ for some $n_1$ between $N$ and $N+1$.  
Set $v:=(p_k,q_k)$.  By the RHS of (\ref{CFA}), the Diophantine 
condition, $|v|\in O(L\log L)$ and $e_0<N\delta$ we get 
$$q_{k+1} \leq |q_k\alpha-p_k|^{-1} \leq (1/c_0)L|v|^{e_0} 
\leq L^{1+N\delta}$$ provided $L\geq L_0$.  
Since $N\leq n_1$, this shows $w$ is normal.  
\hfill\qed

\demo{Proof  of Lemma~{\rm \ref{L:Exist}}}
By the previous lemma, it is enough to prove the existence of 
arbitrarily long $N$-good slits.  We show that a sufficiently long 
slit that is not $N$-good has a nearby slit that is $N$-good.  

Hence, let $w$ be a slit of length $L$ and slope $\alpha$ and 
assume it is not $N$-good.  Let $q_k$ be the largest height in 
$\Spec(\alpha)\cap[1,L^{1+\delta}]$.  Since $q_{k+1}>L^{1+\delta}$ 
(here we use the irrationality of $\lambda$ to guarantee the 
existence of the next convergent) the RHS of (\ref{CFA}) implies 
$\Delta:= (L|q_k\alpha-p_k|)^{-1} > L^\delta$.  
With $L':=L+2mq_k$, it is not hard to see that 
there exists a positive integer $m$ satisfying 
$$e^{N\delta}\log(L')+1/2 \leq \Delta \leq (L')^\delta.$$ Indeed, 
if $m$ is smallest for the RHS, then the LHS holds when $L\geq L_0$.  

Let $w'=w+2mv$ where $v=(p_k,q_k)$.  We show $w'$ is $N$-good.  
Let $\alpha'$ be its slope.  Using $|w'\times v|=|w\times v|$ 
and the cross product formula, we find $|\alpha'-p_k/q_k| = 
1/L'q_k\Delta \leq 1/2q_k^2$ which by (\ref{SCC}) implies 
$q_k\in\Spec(\alpha')$.  Using the above inequalities on $\Delta$ 
in parallel with those in (\ref{CFA}) we obtain 
\begin{eqnarray*}
q_{k+1} &\leq &L'\Delta \leq (L')^{1+\delta} \quad \hbox{and} \\ 
q_{k+1} &\geq& L'(\Delta-q_k/L') \geq e^{N\delta}L'\log L' 
\end{eqnarray*}
which show that $w'$ is $N$-good.  
\enddemo 
 
3.2. {\it  Normal slits have enough normal children.}
Assume $w$ is a normal slit of length $L\geq L_0$ and slope $\alpha$.  
Let $q_k$ be the largest in $\Spec(\alpha)\cap[1,L^{1+\delta}]$ 
and define $n_1\geq1$ uniquely by $q_k=e^{n_1\delta}L\log L$.

\proclaimtitle{Enough children}
\proclaim{Lemma}\label{L:Enough}
Since $w$ has at least $O(L^\delta/\log L)$ $(n-1)$\/{\rm -}\/good children 
where $n:=\min(n_1,N+1)${\rm ,} if $w'$ is a child with length $L'$ 
and slope $\alpha'${\rm ,} then $w'=w+2(p_{k'},q_{k'})$ and 
$q_{k'+1}\in[L'\log L',(L')^{1+\delta}]$ for some $q_{k'}\in\Spec(\alpha')$.  
\endproclaim

{\elevensc Lemma 3.7} (Most children are normal).
{\it The number of children constructed in the previous lemma that 
are not normal is at most \pagebreak $O(L^{\delta-\delta^2}\log L)$.  
}

{\it Proof of Proposition~{\rm \ref{P:Normal}} assuming the 
above lemmas}. A slit $w\in V''$ has enough normal children.  
These are all longer than $L_0$ and therefore  lie in $V''$.  \phantom{socold}\hfill\qed
 
\demo{Proof  of Lemma~{\rm \ref{L:Enough}}} 
Applying Theorem~\ref{Thm3} with $Q=e^{n\delta}L\log L$ and\break
$R=L^{1+\delta}$ gives $O(L^\delta/\log L)$ children of the 
form $w'=w+2(p,q)$ where $\gcd(p,q)=1$, $q\in[L^{1+\delta},
2L^{1+\delta}]$ and $|q\alpha-p|^{-1}\geq e^{n\delta}L\log L$.  

Observe that $L'=L+2q$ so that 
$|\alpha'-p/q| = L|q\alpha-p|/L'q \leq 1/2q^2$.  
Since $\gcd(p,q)=1$, (\ref{SCC}) implies $q=q_{k'}\in\Spec(\alpha')$ 
for some index $k'$.  

It remains to bound the next height $q_{k'+1}\in\Spec(\alpha')$.  
Using the LHS of (\ref{CFA}) together with the lower bound on 
$|q\alpha-p|^{-1}$, we have $$q_{k'+1} \geq L'/L|q\alpha-p|-q_{k'} 
\geq L'(e^{n\delta}\log L - 1/2) \geq e^{(n-1)\delta}L'\log L'$$  
since $L\geq L_0$.  Using the RHS of (\ref{CFA}) and $L\geq L_0$ again 
$$q_{k'+1} \leq L'/Lq|\alpha-p/q| \leq L'L^\delta \leq (L')^{1+\delta}.$$
(In the second step we used the fact that  $|\alpha-p/q|\geq1/L^{2+2\delta}$ 
which, by Lemma~\ref{L:Gaps}, holds for all the children with 
a finite number of exceptions.)  
\enddemo 

\demo{{P}roof of Lemma~{\rm 3.7}} 
Let $\tilde{V}$ be the collection of slits formed by the children 
constructed in Lemma~\ref{L:Enough} that are not normal.  We show 
that $\tilde{V}$ has at most $O(L^{\delta-\delta^2}\log L)$ elements.  
Observe that if $n_1\geq N+1$, all children constructed are $N$-good, 
hence normal, by Lemma~\ref{L:Good=>Normal}.  In this case, $\tilde{V}$ 
is empty and we have nothing to prove.  Therefore, we may assume $n_1<N+1$.  
\vglue9pt 

The next lemma will allow us to count the number of elements in $\tilde{V}$.  
\proclaim{Lemma}\label{L:Reg}%%%%%%%%%%%%%%%%%%%%%%%%%%%%%%%%%%%%%%%%
Let $w'$ be a slit in $\tilde{V}$ of length $L'$ and slope $\alpha'$.  
Then {\rm (i)} the largest $q_{l'}\in\Spec(\alpha')\cap[1,(L')^{1+\delta}]$ 
lies in $[L'\log L',e^{N\delta}L'\log L']$ and {\rm (ii)} it satisfies 
the inequality below for at most finitely many possible values of $a$.  
\vglue12pt
\noindent {\rm (4)} \hfill ${\displaystyle
	\left|L(q_{l'}\alpha-p_{l'}) \pm 2a\right| 
		\leq \frac{1}{L^{n_2\delta+n_2\delta^2}} 
	\quad \hbox{where $n_2:=\max(1,n_1-1)$.}  
}$\hfill
\endproclaim
\vglue8pt

{\it Proof}.
By definition $q_{l'+1}=(L')^{1+n'\delta}$ for some $n'\geq1$.  
Since $w'$ is $(n_1-1)$-good, $q_{l'}=e^{n''\delta}L'\log L'$ 
for some $n''\geq n_1-1\geq0$.  By Lemma~\ref{L:Good=>Normal}, 
$w'$ is not $N$-good so that $n''<N$; this proves (i).  

Since $n_1<N+1$, the fact that $w'$ is not normal implies 
$n'\geq n''\geq n_2$.  

Let $q_{k'}\in\Spec(\alpha')$ be as in Lemma~\ref{L:Enough} 
and recall $q_{k'+1}\in[L'\log L',(L')^{1+\delta}]$.  
By definition of $q_{l'}$, $q_{l'}\geq q_{k'+1}$.  
The recurrence relations satisfied by convergents imply that 
$q_{l'}=aq_{k'+1}+bq_{k'}$ for some integers $a>0$ and $b\geq0$.  
By (i), $a\leq e^{N\delta}$.  

Write $w'=w+2v$ and note that $|w'\times v|=|w\times v|$.  
Since the cross product of consecutive convergents is $\pm1$ 
(thought of   as vectors), 
$$L'|q_{l'}\alpha'-p_{l'}| = |L(q_{l'}\alpha-p_{l'}) \pm 2a|.$$ 
The RHS of (\ref{CFA}), $L'=L+2q\geq L^{1+\delta}$, and $n'\geq n_2$ 
imply 
$$L'|q_{l'}\alpha'-p_{l'}| \leq \frac{L'}{q_{l'+1}} \leq 
	\frac{1}{(L^{1+\delta})^{n'\delta}} \leq 
	\frac{1}{L^{n_2\delta+n_2\delta^2}}$$ 
and (ii) follows.  
\hfill\qed\vglue12pt%%%%%%%%%%%%%%%%%%%%%%%%%%%%%%%%%%%%%%%%%%%%%%%%%%%%%%%

Lemma~\ref{L:Reg} allows us to write $\tilde{V}$ as a finite union of 
subsets $\tilde{V}_{\pm a}$.  Let $Q_{\pm a}$ denote the corresponding 
set of heights $q_{l'}$ associated to the slits in $\tilde{V}_{\pm a}$.  
The next two lemmas complete the proof of Lemma~3.7.
\hfill\qed

\proclaim{Lemma}
$\tilde{V}_{\pm a}$ and $Q_{\pm a}$ have the same number of elements.
\endproclaim

{\it Proof}.
We need to show that the map $\tilde{V}_{\pm a}\rightarrow Q_{\pm a}$ 
sending $w'$ to $q_{l'}$ is injective.  Let $w''$ be different from $w'$ 
with corresponding image $q_{l''}$.  
Note that since $|\alpha'-p_{l'}/q_{l'}|\leq1/(L')^{1+\delta}q_{l'}$ 
is small compared to the distance between the slopes of $w'$ and $w''$ 
(Lemma~\ref{L:Gaps}), the rationals $p_{l'}/q_{l'}$ and $p_{l''}/q_{l''}$ 
are distinct.  Their heights differ because the interval containing 
them is smaller than $1/q_{l'}$.  \phantom{lunchtime}
\hfill\qed

\proclaim{Lemma}
Each $Q_{\pm a}$ is a union of at most $O(\log L)$ subsets{\rm ,} 
each having at most $O(L^{\delta-\delta^2})$ elements.  
\endproclaim

{\it Proof}.
Let $q_{l'},q_{l''}\in Q_{\pm a}$ and set $\bar{q}:=|q_{l''}-q_{l'}|$.  
We claim that if $\bar{q}\leq L^{1+\delta}$ then $\bar{q}=dq_k$ 
for some positive integer $d\in O(L^{\delta-\delta^2})$.  
This implies the lemma since we either have $\bar{q}>L^{1+\delta}$ 
or $\bar{q}\ll L^{1+\delta}$ so that the elements of $Q_{\pm a}$ fall 
into $O(\log L)$ clusters, each having $O(L^{\delta-\delta^2})$ elements.  

Hence, assume $\bar{q}\leq L^{1+\delta}$ and set $\bar{p}=|p_l-p_{l'}|$.  
Let $p/q$ be the reduced form of $\bar{p}/\bar{q}$ 
so that $\bar{q}=dq$ where $d=\gcd(\bar{p},\bar{q})$.  
From (4) and the triangle inequality 
$$|\bar{q}\alpha-\bar{p}| \leq \frac{2}{L^{1+n_2\delta+n_2\delta^2}}.$$
Since $q\leq\bar{q}\leq L^{1+\delta}$, $|\alpha-p/q|\leq1/2q^2$.  
By (\ref{SCC}) $p/q$ is a convergent of $\alpha$; 
since $\gcd(p,q)=1$, $q=q_{k'}\in\Spec(\alpha)$ for some index $k'$.  
By hypothesis, $q\leq L^{1+\delta}$, so we must have $k'\leq k$.  
In fact, we must have $k'=k$ because $k'<k$ implies $|\bar{q}\alpha-\bar{p}|
=d|q_{k'}\alpha-p_{k'}| \geq 1/(q_k+q_{k-1}) \geq 1/2q_k$ which 
contradicts the previous inequality.  Using the LHS of (\ref{CFA}), 
we now have 
$$d = \frac{|\bar{q}\alpha-\bar{p}|}{|q_k\alpha-p_k|} \leq 
      \frac{2(q_{k+1}+q_k)}{L^{1+n_2\delta+n_2\delta^2}}$$ 
which is $O(L^{(n_3-n_2)\delta-n_2\delta^2})$ 
where $n_3$ is defined by $q_{k+1}=L^{1+n_3\delta}$.  
Since $w$ is normal, $n_3\leq n_1$.  This together with (4) 
implies $n_3-n_2\leq n_1-n_2\leq1$ and $n_2\geq1$; thus proving the claim.  
\hfill\qed

\section{Counting rationals in intervals}\label{S:Density}

Primitive lattice points in the parallelogram $\Sigma$ correspond to 
rationals in 
the interval $I:=\left[\alpha-\frac{1}{RQ},\alpha+\frac{1}{RQ}\right]$.    
Set $$\Lambda_I := \left\{ (x,y)\in{\Bbb R  }^2 : 
		x/y\in I,\; R \leq y \leq 2R \right\}.$$ 

\specialnumber{4}\proclaim{Theorem} \label{Thm4}
If $\Spec(\alpha)\cap[Q,R]\neq\emptyset${\rm ,} $\dens(\Lambda_I)\geq1/24$ 
and $R/Q\geq 16$.  
\endproclaim
\advance\eqcount by 4

{\it Proof of Theorem~{\rm \ref{Thm3}}}.
We use $A_0=16$ and $\rho_0=1/32$.  Observe that $\Sigma$ contains 
$\Lambda_{I'}$ where $I'$ is concentric with $I$ and half as wide.  
Moreover, $\Lambda_{I'}$ occupies three quarters of its area.  
Replacing $Q$ with $2Q$ in Theorem~\ref{Thm4} and assuming 
$\Spec(\alpha)\cap[2Q,R]\neq\emptyset$, we conclude 
$\dens(\Sigma)\geq(3/4)\dens(\Lambda_{I'})\geq\rho_0$.  

If $\Spec(\alpha)\cap[2Q,R]=\emptyset$, let $q_k$ be the largest 
height in $\Spec(\alpha)\cap[Q,R]$.  By the RHS of (\ref{CFA}) 
it is easy to show the reduced fraction 
	$$\frac{p}{q} = \frac{ap_k+p_{k-1}}{aq_k+q_{k-1}} 
		\quad \hbox{where} \quad a=1,2,3,...$$ 
(known as an intermediate fraction of $\alpha$ when $a\leq a_{k+1}$) 
satisfies 
\begin{equation}\label{E:Ifrac}
    \left|\alpha-\frac{p}{q}\right| \leq 
	\left(\frac{q_k+|q_{k+1}-q|}{q}\right)\frac{1}{q_kq_{k+1}}.  
\end{equation}
These correspond to $R/q_k$ primitive lattice points in $\Sigma$, and 
since $q_k\leq2Q$, it follows easily that $\dens(\Sigma)\geq1/4\geq\rho_0$.  
\hfill\qed

\proclaim{Lemma}\label{L:Rint}  
 If $J$ has 
rational endpoints of height at most $R${\rm ,} $\dens(\Lambda_J)\break\geq1/6$.  
\endproclaim

{\it Proof}.
We shall first prove the lemma under the additional hypotheses:
\begin{itemize}
    \item[(i)] the height of any rational in ${\rm int}(J)$ 
		is greater than $R$, and 
    \item[(ii)] $|pq'-p'q| = 1$, where $p/q$ and $p'/q'$ are 
		the endpoints of $J$.  
\end{itemize}
By (ii), arithmetic density is preserved by the linear map $\gamma$ 
which sends the standard basis to lattice points corresponding to 
the endpoints of $J$.  
Note that $$\gamma^{-1}(\Lambda_J)\subset\Delta:=
\left\{ (x,y)\in{\Bbb R  }^2 : x/a+y/a'\leq2, \; x,y>0 \right\}$$ 
where $a=R/q$ and $a'=R/q'$.  
Let $n(\Delta)$ denote the number of primitive lattice points 
${\rm int}(\Delta)$.  By (i), it is enough to show $n(\Delta)\geq aa'/4$.  

Without loss of generality, assume $a'\geq a\geq1$.  There are two cases.  
If $a'\leq2$ then $aa'\leq4$, and since $(1,1)\in\Delta$ we have 
$n(\Delta)\geq1\geq aa'/4$.  On the other hand, if $a'>2$, then since 
$\gamma(1,1)\in\Lambda_J$, we have $1/a+1/a'\geq1$ so that $a\in[1,2]$.  
Considering pairs of the form $(1,n)$, we find $$n(\Delta)
\geq \left\lfloor a'\left(2-\frac{1}{a}\right) \right\rfloor 
> aa'\left(\frac{2}{a}-\frac{1}{a^2}-\frac{1}{2}\right)
\geq \frac{aa'}{4}.$$  
This completes the proof assuming the additional hypotheses.  

Note that every interval is a disjoint union of intervals satisfying (i).  
Hence, the lemma follows if we show ${\rm (i)}= {\rm (ii)}$.  Indeed, let $d=|pq'-p'q|$.  
There is a linear map in ${\rm GL}_2{\Bbb Z}$ that takes $(p,q)$ to $(0,1)$ 
and $(p',q')$ to $(d,d')$ for some integer $d',0\leq d'<d$.  If $d'>0$ 
then $(1,1)$ is contained in the triangle determined by the origin, 
$(1,0)$ and $(d,d')$ and corresponds to a rational of height at most $R$ 
in ${\rm int}(J)$.  Therefore, (i) implies $d'=0$, and since $\gcd(d,d')=1$, 
this in turn implies that $d=1$, giving (ii).  
\hfill\qed\vglue12pt

The height of a rational strictly between $p/q$ 
and $p'/q'$ is at least $q+q'$: 
$$\frac{1}{qq'} = \left|\frac{p}{q}-\frac{p''}{q''}\right| 
    + \left|\frac{p''}{q''}-\frac{p'}{q'}\right| 
    \geq \frac{1}{qq''} + \frac{1}{q'q''}$$
This will be used several times in the next proof.  

\demo{Proof of Theorem~{\rm \ref{Thm4}}}
Let $\alpha'$ and $\alpha''$ denote the left and right endpoints 
of $I$, respectively.  Let $q_k$ be the largest height in 
$\Spec(\alpha)\cap[Q,R]$.  By (\ref{CFA}), $p_k/q_k\in I$ and 
without loss of generality, we assume $p_k/q_k\leq\alpha$.  
Let $q_l$ be the largest height in $\Spec(\alpha'')\cap[1,R]$.  
Since $p_k/q_k$ cannot lie strictly between $p_l/q_l$ and 
$p_{l+1}/q_{l+1}$, it must be the case that $p_k/q_k\leq p_l/q_l$.  
In fact, strict inequality must hold because $q_l=q_k\geq Q$ and 
$|\alpha'-p_k/q_k|\leq1/q_lq_{l+1}<1/RQ$ give a contradiction.  

Let $J:=[p_k/q_k,p_l/q_l]$.  We claim its length is at least $1/2RQ$.  
In fact, $p_l/q_l$ lies within $1/2RQ$ of $\alpha''$ if $q_l\geq2Q$.  
On the other hand, $q_l\leq2Q$ implies $|J|\geq1/q_kq_l\geq1/2RQ$.  
In either case, $|J|\geq1/2RQ$.  

We may assume $\alpha''\leq p_l/q_l$, for otherwise $J\subset I$ 
and $\dens(\Lambda_I)\geq1/24$.  

There are three cases.  
First, if $q_{l+1}+3q_l>2R$ there can be at most two rationals 
with height at most $2R$ that lie strictly between $p_l/q_l$ 
and $p_{l+1}/q_{l+1}$:
$$\frac{p_{l+1}}{q_{l+1}} < \frac{p_{l+1}+p_l}{q_{l+1}+q_l} 
  < \frac{p_{l+1}+2p_l}{q_{l+1}+2q_l} < \frac{p_l}{q_l}.$$  
Let $n(\Lambda_I)$ denote the number of primitive lattice 
points in $\Lambda_I$.  By Lemma~\ref{L:Rint}, 
$$n(\Lambda_I) \geq \frac{\area(\Lambda_J)}{6}-2= 
  \left(\frac{1}{12}-\frac{2Q}{3R}\right)\area(\Lambda_I)$$
and since $R/Q\geq16$, $\dens(\Lambda_I)\geq1/24$.  

Next, suppose $q_{l+1}+3q_l\leq2R$ and $q_l\geq2Q$.  
Note that an intermediate fraction of $\alpha''$ with height 
between $q_l$ and $q_{l+1}$ lies to the left of $\alpha''$.  
Let $J':=[p_k/q_k,p/q]$ where $p/q$ is the intermediate fraction 
with the largest height not exceeding $R$.  By definition, $q>R-q_l$.  
From (\ref{E:Ifrac}) we have 
$$\left|\alpha-\frac{p}{q}\right| 
	\leq \left(\frac{q_l+q_{l+1}-q}{q}\right)\frac{1}{q_lq_{l+1}} 
	\leq \left(\frac{2R-2q_l-q}{R-q_l}\right)\frac{1}{2RQ} 
	\leq \frac{1}{2RQ}$$ 
so that $|J'|\geq1/2RQ$ and $\dens(\Lambda_I)\geq1/24$.  

Finally, assume that $q_{l+1}+3q_l\leq2R$ and $q_l<2Q$.  
Again, we consider the intermediate fractions of $\alpha''$.  
Observe that their heights are at most $2R$, since $q_{l+1}\leq2R$.  
They form a sequence that increase towards $\alpha''$ from the left.  
Given a consecutive pair with heights less than $R$ we can always 
find a rational strictly in between them with height in $[R,2R]$.  
It follows that the number of rationals in $I$ with height in $[R,2R]$ 
is at least $(q_{l+1}-q)/q_l$, where $q$ is the height of the first 
intermediate fraction that falls into $I$.  According to (\ref{E:Ifrac}), 
an intermediate fraction lies in $I$ as soon as its height is greater than 
$$R' := \frac{(q_l+q_{l+1})RQ}{2q_lq_{l+1}+RQ}.$$  
Hence, $q-q_l\leq R'$ and 
\begin{eqnarray*}
n(\Lambda_I)&\geq&\frac{q_{l+1}-q}{q_l}\geq\frac{q_{l+1}-R'}{q_l}-1 
=\frac{2q_{l+1}^2-RQ}{2q_lq_{l+1}+RQ}-1\\ &\geq& \frac{2R^2-RQ}{9RQ}-1 
=\frac{2}{27}\left(1-\frac{5Q}{R}\right)\area(\Lambda_I).  
\end{eqnarray*}
Since $R/Q\geq16$, we get $\dens(\Lambda_I)\geq1/24$.  
\enddemo
  \vglue12pt
\vbox{\baselineskip 12pt \eightsc Northwestern University, Evanston, IL 

{\eightpoint {\it E-mail address\/}: yitwah@math.northwestern.edu}}

 \def\bR{{\Bbb R}}
\def\bC{{\Bbb C}}
\def\bN{{\Bbb N}}
\def\bP{{\Bbb P}}
\def\bQ{{\Bbb Q}}
 
\def\cB{{\cal  B}}
\def\cL{{\cal  L}}
\def\cM{{\cal  M}}
\def\cB{{\cal  B}}
\def\cC{{\cal  C}}
\def\cO{{\cal  O}}
\def\cF{{\cal  F}}
\def\cG{{\cal  G}}
\def\cH{{\cal  H}}
\def\cE{{\cal  E}}
\def\cU{{\cal  U}}
\def\a{\alpha}
\def\l{\lambda}
\def\L{\Lambda}
\def\fs{F_{\sigma}}
\def\gd{G_{\delta}}
\def\bsq{$\blacksquare$}
\def\qqquad{\hspace{30mm}}
\def\theend{\\[-7mm] \begin{flushright}
  \bsq \hspace{7mm} \  \end{flushright}}
\newcommand{\az}[1]{\mbox{{ \bsq\,{\bf #1}\,\bsq}}}
\def\defeq{\stackrel{{\rm def}}{=}} %%% define by def

\vglue-.5in

\title{Appendix}
 \shorttitle{Appendix}
\vglue-12pt
   \author{Michael Boshernitzan}

With notation as in the introduction, for  
  $\l\in [0,1)$,  set   
\begin{equation}\label{E:h} 
   h(\l) \defeq \hbox{H.dim NE}(Q_\l).
\end{equation}           

 Recall (see Introduction) that 
   $h(\l)\!\le\!1/2$ 
   for all $\l$ \cite{Ma},
   and, by Theorem~\ref{Thm1}, that $h(\l)=1/2$
   for all  Diophantine  $\l$.
   We also have  $h(\l)=0$  for rational  $\l$
    (then the set  $Q_\l$  is in fact countable
    \cite{V1}). The main result in this section is given by the 
    following theorem.
 
\specialnumber{5}\proclaim{Theorem}  \label{T:MB0}
  The set of $\l\in [0,1)$  for which 
  ${\rm H.dim\ NE}(Q_\l)=0$
  form a residual subset of  $[0,1).$
  In particular{\rm ,} there are irrational  
  $\l\in [0,1)$  such that  $h(\l)=0$.  
\endproclaim
 
Recall that a subset  $A \subset X$  of 
  a topological space  $X$  is called {\em residual} 
  (or topologically  large)  if it contains a dense 
  $\gd$-subset of  $X$. A subset $Y\subset X$  is called
  a $\gd$-set (in  $X$)  if $Y$  is a countable
  intersection of open subsets of  $X$. Its complement 
  $X\setminus Y$  is called an $\fs$-set. 
  
We remark that no irrational  number $\l$
  satisfying  $h(\l)=0$
  is known even though  the set is topologically large 
  (in particular, uncountable). Any such  $\l$  must be Liou\-ville.  
   (Note that the set of Liouville numbers  forms
   a residual set of Lebesgue measure~0.)

\vglue16pt
\centerline{\bf A.1.   ${\Bbb Z}_2$ skew products of irrational rotations}
\vglue12pt

Let  $X= S^1_0 \cup S^1_1= S^1 \times  \{0,1\}$ be the union 
  of two unit circles  $S^1_k=S^1 \times \{n\}$, \mbox{$n\in\{0,1\}$},
  and consider the two-parameter family    of transformations 
$$ %beg
    \rho _{\a,\l}: X\to X,  \hspace{12mm}
    \a\in \bR,\; \l\in K=[0,1),
$$  
defined as follows.  For  $x=(s,n)\in X$, 
\begin{equation} \label{e:ro}
     \rho_{\a,\l}(x)=\rho_{\a,\l}(s,n)=(s\oplus\a, n'),  \qquad  
     n'=\left\{
  \begin{array}{ll}
      \!\! n,&  \hbox{if }\;   0 \le s < \l, \\
      \!\! 1-n,&  \hbox{if }\;  \l \le s <1,
    \end{array}
    \right.
\end{equation}
where 
   $x=(s,n) \in X, \quad 
     s,s\oplus \a \in S^1=\bR /{\Bbb Z} = [0,1),$ 
and  $\oplus$ stands for the group operation in  $S^1$. 

The dynamical systems  $(X,\rho _{\a,\l})$  have been studied
by Veech [V1] as particular  ${\Bbb Z}_2$  skew product  extensions 
of irrational $\a$-rotations.  Indeed, $\rho_{\a,\l}$ may be 
interpreted as the first return map to a disjoint union of two 
circles embedded in the surface associated to $Q_\lambda$.  
In particular, properties of billiards on  $Q_\l$  reduce  to 
the study of dynamical systems  $(X,\rho _{\a,\l})$.  
One verifies that the ergodicity of the billiard system $Q_{\l}$ 
in direction $\theta$ is equivalent to the ergodicity of the map 
$\rho_{\a,\l}$ with $\a=\tan (\theta)$ (the slope in direction 
$\theta$).  Denote 
\def\rne{\left(X,\rho _{\a,\l}\right)  \hbox{ is not ergodic}}
$$
    \hbox{NE}(X) = \left\{ (\a,\l) \in \bR\times [0,1)\mid \rne \right\},
$$
and,  for $\l \in [0,1)$, 
\begin{eqnarray} \label{E:NEX}
    \hbox{NE}(X_\l) &=&  \left\{ \a \in \bR\mid (\a,\l) \in \hbox{NE}(X)
         \right\}   \\
	& =&  \left\{ \a \in \bR\mid  \rne  \right\}.
\nonumber\end{eqnarray}

The sets \  NE$(X_\l)$, NE$(Q_\l)$
have the same  Hausdorff dimensions because 
$$\hbox{NE}(X_\l)=\hbox{NE}(X_\l)+1=\tan(\hbox{NE}(Q_\l))=
     \left\{ \tan(\theta) \mid \theta\in \hbox{NE}(Q_\l)\right\}.
$$
Thus we have  (see \eqref{E:h})
\begin{equation}   \label{E:hl}
     h(\l)=\hbox{H.dim NE}(Q_\l)=\hbox{H.dim NE}(X_\l).
\end{equation}

\centerline{\bf A.2.   The  topological lemma}
\vglue12pt

  The following lemma is central in the proof of
Theorem \ref{T:MB0}. 

\specialnumber{A.1}\proclaim{Lemma} \label{L:Top}
    Let  $L$  be a $\gd$\/{\rm -}\/subset of a $\sigma$\/{\rm -}\/compact metric
    space~$K$.  Let  $P$   be a Polish space
    {\rm (}\/a space with a complete metric topology\/{\rm ).}  Let  $H$  be an
    $\fs$\/{\rm -}\/subset of the cartesian product  $W=K \times P$.  
    For every  $p \in P${\rm ,}  denote
\begin{equation}  \label{E:KP}
     K(p) = \{ k \in K \mid (k,p) \in H \} \subset K
\end{equation}
 and  
\begin{equation}  \label{E:PO}
     P_{o} = \{ p \in P \mid K(p) \subset L \}.
\end{equation}
If $P_{o}$   is dense in  $P${\rm ,}  then  $P_{o}$  is a 
residual subset {\rm (}i.e.{\rm ,} contains a dense\break  $\gd$\/{\rm -}\/subset\/{\rm )}
of $P$. In particular{\rm ,}   $P_{o}$  is
uncountable if  $P$  is.
\endproclaim

{\it Proof}.
Since the family of residual  subsets of  $P$ is closed 
  under coun\-table intersections, we assume 
  (as we may without loss of generality)  that
  $K$  is compact and $H$  is closed 
  in $W=K \times P$.  

  Since  $P$  is separable and  $P_{o}$ is dense in  $P$,  
there is a countable subset  $P_{c} \subset P_{o}$ which 
is dense in $P$. We assume that the points 
of   $P_c$  are arranged in one  sequence  
$\left\{p_i\right\}$   so that every 
point  $p \in P_{c}$  is repeated infinitely many times, 
i.e.   $p_{i}=p$, for infinitely many $i \ge 1$. 

  Let  $\pi_P\!\!: W\to P$, $\pi_K\!\!: W\to K$  be
canonical projections.  Denote 
\begin{eqnarray}
K_i&=&K(p_i)=\pi_{K}(\pi^{-1}_{P}(p_i) \cap H),
\\[6pt]
K'_i&=&K_i \times \{p_i\} = \pi^{-1}_{P}(p_i)\cap H) \subset 
    L \times \{p_i\}.
\end{eqnarray}

  Since  $L$  is a $\gd$-subset of $K$,  $L$ has a representation 

$$
L=\bigcap_{i \ge 1} U_i, \qquad \;
   K \supset U_1\supset U_2 \supset U_3 \ldots
$$
\noindent  where $\{ U_i \}$ forms  
(without loss of generality) a nonincreasing sequence of open 
subsets of  $K$.  

  Fix any integer  $i \le 1$.  The set  
  $H_i = H \setminus  \pi_K^{-1}(U_i) \subset W$ is closed,  and so is the set   $\pi_P(H_i) \subset P$  (the 
  projection  $\pi_P \! : W\to P$  is a closed map since   
  $K$ is  compact).  One observes that   
  $p_i\notin \pi_P(H_i´)$
  (since  $K_i´ \subset L \subset U_i´$),  and hence
  $$
      \pi_P^{-1}(p_i) \cap H = K'_{i} =  K_{i} \times \{ p_i \} 
      \subset    (U_{i} \times \{ p_i \} ) \cap H \subset
      \pi_K^{-1} \cap H.
  $$

  Therefore the set   
  $V=\bigcap_{ k \ge 1} ( \bigcup_{ i \ge k} V_i )$
  is a dense  $\gd$-subset of  $P$ (it contains a dense subset \pagebreak  $P_c$).  
  
    To complete the proof,  we have to verify that   $V \subset P_{o}$.
  If  $v \in V$,  then  $v \in V_i$,  for an infinite set of  $i$. 
  For all those  $i$,  
  $$
       v \notin \pi_P(H_i)=\pi_P (H_i \cap 
       \pi_K^{-1}(K \setminus U_i))  
  $$
  which implies   $K(v) \subset U_i$.   It follows that
  $K(v) \subset L = \bigcap_{ i \ge 1} U_i$,  and therefore   
  $v \in  P_{o}$. 
  This completes the proof of Lemma \ref{L:Top}.  

\vglue16pt \centerline{\bf A.3.  The completion of the proof}  
\vglue12pt

 Let  $\cH$  be a separable Hilbert space.  
 Denote by  $\cL(\cH)$ the Banach algebra of bounded 
 linear operators on  $\cH$, and by  $\cU(\cH)$  the 
 subset of unitary operators  on~$\cH$.
 We recall that convergence  $T_i \to T$  in the strong 
   (or weak) operator topology means 
   convergence    $T_if \to Tf$,  for all  $f\in\cH$,
   in the strong (or weak, respectively)  topology of 
   the Hilbert  space  $\cH $.   
   The strong and weak operator topologies on  
   $\cL(\cH)$ are  different,  but they coincide 
   when restricted to the set  \mbox{$\cU(\cH) \subset \cL(\cH)$}
   (see e.g. Halmos~\cite[pp.\ 61--80]{Ha}). 
   
   Now  we view $X= S^1 \times  \{0,1\}$  as
   a probability measure space  $(X,\cB,\mu)$ with  
   $\mu$   a multiple of the Lebesgue measure, $d\mu =ds/2$. 
   Take  $\cH =L^2(X,\cB,\mu)$,  and denote by  $\cG(X)$ the
   family of invertible measure-preserving transformations
   $T: X\to X$.  Then  $\cG(X)$  is naturally imbedded into
   $\cU(\cH)$:  for  $T\in \cG(X)$  and  $f\in\cH=L^2(X,\cB,\mu)$
   define   $T(f(x))=f(T(x))\in\cH$.  The subspace topology 
   on  $\cG(X)$  induced by a strong  (equivalently, weak) operator
   topology on   $\cU(\cH)$  is called  weak topology on  $\cG(X)$. 
   The following result is well known (see \cite[p.\ 80]{Ha}).
   
\specialnumber{A.2} \proclaim{Lemma} \label{L:halmos}
     The set 
   $$
      \cE(X)=\{T \in \cG(X) \mid T \textup{ \ is ergodic} \}\subset\cG(X)
   $$
   is a dense  $\gd$\/{\rm -}\/subset of   $\cG(X)$
   {\rm (}\/in the weak topology of  $\cG(X)$\/{\rm ).}\/
 \endproclaim
   
 Now let  $K=S^1$ and  $P=[0,1)$,  and define 
  $W=K \times P=S^1\times [0,1)$, just as in 
  Lemma \ref{L:Top}.  
  The  map    $\phi: W \to \cG(X)$  defined by the formula
  (see  \eqref{e:ro})
 $$
  \phi(w)=\rho_{w}=\rho_{k,p}, 
  \hspace{6mm} \hbox{for } \ w=(k,p) \in W=K \times P
 $$
   is easily verified to be continuous.
   In view of Lemma \ref{L:halmos}, 
    the set   $\phi^{-1}(\cE(X))$  is a  $\gd$-subset of   
    $W$,  and thus the complement
    $$
    H \defeq W \setminus \phi^{-1}(\cE(X)) =
       \{w=(k,p)\in W \mid  \phi(w)=\rho_{k,p} \textup{ \ is not ergodic}\}
    $$
    is an $\fs$-subset of  $W$. 
    
    Denote by  $\bQ(P)$,  $\bQ(K)$ the sets of rationals 
      in~$P$  and $K=\bR/{\Bbb Z}=[0,1)$, respectively.    Fix an  arbitrary  
       $\gd$-subset  $L$  in  $K$  such that 
\begin{equation}\label{E:L}
               \bQ(K)\subset L\subset  K, \qquad \hbox{H.dim }L=0.
\end{equation} 
    
To satisfy the conditions of  Lemma \ref{L:Top},  it remains to 
   verify that  $P_{o}$  is dense in  $P$.
   For every  $\l\in P$, we have  $K(\l)=$ NE$(X_\l)$  
      (see \eqref{E:NEX}, \eqref{E:KP}), and therefore
$$
   P_{o}=\{ \l\in P \mid \hbox{NE}(X_\l) \subset L\}.
$$ 

For $\l\in\bQ(P)$, Veech \cite{V1} proved that NE$(X_\l)=\bQ(K)$.
    Since $L\supset \bQ(K)$  by the choice of  $L$,  the inclusion
    $\bQ(P)\subset P_{o}$ holds. Thus  $P_{o}$ is indeed dense in~$P$, and, by
    Lemma \ref{L:Top},   $P_{o}$  is a residual subset of  $P=[0,1)$.
    
For every  $\l\in P_{o}$, we have  NE$(X_\l) \subset L$;  thus 
   H.dim NE$(X_\l)=0$  in view of \eqref{E:L}.  This completes the 
   proof of Theorem \ref{T:MB0}.

\AuthorRefNames [KMS]


\begin{references}


\bibitem{Fa}
  \name{K.\ Falconer}, 
   {\it Fractal Geometry\/}. {\it Mathematical Foundations and
   Applications\/}, 
        John Wiley \& Sons Ltd., Chichester, 1990. 

\bibitem{Fu}
  \name{H.\ Furstenberg}, 
  Strict ergodicity and transformation of the torus, 
        {\it Amer.\ J.\ Math\/}.\  {\bf 83} (1961), 573--601.  

\bibitem{Ha}
  \name{P.\ Halmos},
   Lectures on ergodic theory,
  {\it Math.\ Soc.\ of Japan},
    no.\  3 (1956), 1--99.

\bibitem{Kh}
  \name{A.\ Ya.\  Khintchine}, 
  {\it Continued Fractions}, translated by Peter Wynn, 
        P.~Noordhoff Ltd., Groningen, 1963.

\bibitem{KMS} 
  \name{S.\ Kerckhoff, H. Masur}, and \name{J.\ Smillie}, 
   Ergodicity of billiard flows and quadratic differentials, 
        {\it Ann.\ of Math\/}.\ {\bf 124} (1986), 293--311.  

\bibitem{Ma} 
  \name{H.\ Masur}, 
   Hausdorff dimension of the set of nonergodic foliations of a 
        quadratic differential, 
        {\it Duke Math.\ J\/}.\  {\bf 66} (1992), 387--442.  

\bibitem{MS} 
  \name{H.\ Masur} and \name{J.\ Smillie}, 
   Hausdorff dimension of sets of nonergodic measured foliations, 
        {\it Ann.\ of Math\/}.\ {\bf 134}  (1991), 455--543.  

\bibitem{MT}
  \name{H.\ Masur} and \name{S.\ Tabachnikov}, 
  Rational billiards and flat structures, in 
  {\it Handbook of Dynamical Systems\/} 
{\bf 1A\/} (2002), 1015--1089.

\bibitem{V1}
  \name{W.\ Veech}, 
   Strict ergodicity in zero dimensional dynamical systems and 
        the Kronecker-Weyl theorem mod $2$, 
        {\it Trans.\ Amer.\ Math.\ Soc\/}.\  {\bf 140} (1969), 1--34.  

\bibitem{V2}
  \bibline,
   The billiard in a regular polygon, 
	{\it Geom.\ Funct.\ Anal\/}.\  {\bf 2} (1992), 341--379.  

\end{references}
\end{document}